# Analytical solution with tanh-method and fractional sub-equation method for non-linear partial differential equations and corresponding fractional differential equation composed with Jumarie fractional derivative

Uttam Ghosh (1), Srijan Sengupta (2a), Susmita Sarkar (2b), Shantanu Das (3)

(1): Department of Mathematics, Nabadwip Vidyasagar College, Nabadwip, Nadia, West Bengal, India;
Email: uttam_math@yahoo.co.in
(2): Department of Applied Mathematics, University of Calcutta, Kolkata, India
Email (2b): susmita62@yahoo.co.in,
(3)Scientist H+, RCSDS, Reactor Control Division BARC Mumbai India
Senior Research Professor, Dept. of Physics, Jadavpur University Kolkata
Adjunct Professor. DIAT-Pune
UGC-Visiting Fellow. Dept of Appl. Mathematics; Univ. of Calcutta, Kolkata
Email (3): shantanu@barc.gov.in

**Abstract**

The solution of non-linear differential equation, non-linear partial differential equation and non-linear fractional differential equation is current research in Applied Science. Here tanh-method and Fractional Sub-Equation methods are used to solve three non-linear differential equations and the corresponding fractional differential equation. The fractional differential equations here are composed with Jumarie fractional derivative. Both the solution is obtained in analytical traveling wave solution form. We have not come across solutions of these equations reported anywhere earlier.

## 1.0 Introduction

The advanced analytical methods to get solution to non-linear partial differential equation is current research field in different branches of applied science such as Applied Mathematics, Physics, Biology and Engineering applications. On the other hand fractional calculus and consequently solutions of fractional differential equations is also another growing field of science especially in Mathematics, Physics, Engineering and other scientific fields like controls etc [1-5, 25]. The exact solutions of non-linear partial differential equation and the corresponding fractional differential equations are physically very important [25-27]. To solve such type non-linear differential equations different methods are used, such as Adomain Decomposition Method (ADM) [6-7, 25], the Variational Iterative Method (VIM) [8-10], the Homotopy Perturbation Method (HPM) [11-13], Differential Transform Method (DTM) [14]. Another important method to find the traveling wave solution of the non-linear partial differential equation is tanh- method was first introduced by Huibin and Kelin [15]. This method was used directly into a higher-order KdV equation, to get its solutions. Wazwaz [16, 17] used tanh-method to find traveling wave solutions of some non-linear equations generalized Fisher's equation, generalized KdV equation, Zhiber-Shaba equation and other related equations. Recently Zhang and Zhang [18] introduce Fractional Sub-Equation Method to find travelling wave solution of the non-linear fractional differential equation. The principle of this method is expressing the solution of the non-linear differential equation as polynomials, where the variables of which is one of the solution of simple and solvable ordinary and partial differential equation. The method is also based on Homogeneous Balance Principle [19, 20] and applied on Jumarie's modified Riemann-Liouvelli derivative [21, 22]. The Fractional Sub-Equation Method used by many authors [19, 21, 24] to find the analytic solution of non-linear fractional differential equations. Rest of the paper is as follows-in section 2.0 we describe tanh-method, in section 3.0 we describe Fractional Sub-Equation method, in section 7.0-12.0 we describe the solution of three non-linear equations in both integer order partial differential and fractional order partial differential equation types. In this paper we find the solution of three non-linear equation using tanh-method and the corresponding non-linear fractional differential equation using Fractional Sub-Equation method. In all the cases analytical solutions obtained in travelling wave solution form. In this paper the symbols for fractional differential operator used $_0^J D_x^\alpha$, $f^{(\alpha)}$ means Jumarie fractional derivative.

## 2.0 Tanh-method

Let the non-linear partial differential equation is of the form

$$H(u, u_t, u_x, u_{xx}....) = 0 ..............................................................(1)$$

Where $u = u(x,t)$ is function of $x$, $t$.

For instant

$$\frac{\partial^2 u}{\partial x^2} - u \frac{\partial u}{\partial t} = 0$$

is of such an equation as $H(u, u_t, u_{xx}) = 0$, with $H$ a linear /non-linear operator. Using travelling wave transformation $\xi = kx + ct$, with $k$ constant called wave-number (or propagation constant) and $c$ a constant; as velocity of propagating wave, equation (1) reduce to an ordinary differential equation in the form

$$H(u, u_\xi, u_{\xi\xi}....) = 0 ..............................................................(2)$$

Where $u = u(\xi)$, is function of one variable $\xi$. In summary we are making several variables $(x, y, z, t)$ to a single variable $\xi$, with travelling wave transformation such as $\xi = k_x x + k_y y + k_z z + ct$, with corresponding constants $k_x$, $k_y$, $k_z$, $c$ as described above, and $\xi$ is in units of distance. Thus the operator $\frac{\partial}{\partial t}$, in terms of $\xi$ with $\xi = kx + ct$ is following

$$\frac{\partial}{\partial t} \equiv \frac{\partial}{\partial \xi}\frac{\partial \xi}{\partial t} = c\frac{d}{d\xi}$$

and similarly $\frac{\partial^2}{\partial x^2}$ is

$$\frac{\partial^2}{\partial x^2} \equiv \frac{\partial}{\partial x}\left(\frac{\partial}{\partial x}\right) = \frac{\partial}{\partial \xi}\frac{\partial \xi}{\partial x}\left(\frac{\partial}{\partial \xi}\frac{\partial \xi}{\partial x}\right) = \frac{\partial}{\partial \xi}k\left(\frac{\partial}{\partial \xi}k\right) = k^2\frac{d^2}{d\xi^2}$$

The transformed equation $\frac{\partial^2 u}{\partial x^2} - u\frac{\partial u}{\partial t} = 0$, is thus,

$$k^2 \frac{d^2 u}{d\xi^2} - cu\frac{du}{d\xi} = 0 \quad \text{or} \quad k^2 u_{\xi\xi} - cuu_\xi = 0$$

which is $H(u, u_\xi, u_{\xi\xi}) = 0$.

Again, let $\phi = \phi(\xi)$ be such type of function which is a solution of Riccati equation $\phi'(\xi) = \sigma + \phi^2$. Where $\sigma$ is a real constant. That is $\frac{d\phi}{d\xi} = \sigma + \phi^2$, $\int d\xi = \int \frac{d\phi}{\sigma + \phi^2}$

Then

$$\phi(\xi) = \begin{cases} \left.\begin{array}{l}-\sqrt{-\sigma}\tanh(\sqrt{-\sigma}\xi) \\ -\sqrt{-\sigma}\coth(\sqrt{-\sigma}\xi)\end{array}\right\} & \text{for } \sigma < 0 \\ \left.\begin{array}{l}\sqrt{\sigma}\tan(\sqrt{\sigma}\xi) \\ -\sqrt{\sigma}\cot(\sqrt{\sigma}\xi)\end{array}\right\} & \text{for } \sigma > 0 \\ -\dfrac{1}{\xi} & \text{for } \sigma = 0 \end{cases}$$

The above expression is detailed below briefly.

Case-I
For $\sigma < 0$ then the Riccati equation can be written equation, considering integral constant as zero

$$\int d\xi = \int \frac{d\phi}{\sigma + \phi^2} \qquad \xi = \begin{cases} -\frac{1}{\sqrt{-\sigma}} \tanh^{-1}\left(\frac{\phi}{\sqrt{-\sigma}}\right) \\ -\frac{1}{\sqrt{-\sigma}} \coth^{-1}\left(\frac{\phi}{\sqrt{-\sigma}}\right) \end{cases}$$

Therefore, by inverting above we get the following

$$\phi = \begin{cases} -\sqrt{-\sigma} \tanh\left(\sqrt{-\sigma}\,\xi\right) \\ -\sqrt{-\sigma} \coth\left(\sqrt{-\sigma}\,\xi\right) \end{cases}$$

Case-II

For $\sigma = 0$ then the Riccati equation can be written equation, ignoring integration constant

$$\int d\xi = \int \frac{d\phi}{\phi^2} \qquad \xi = -\frac{1}{\phi} \qquad \phi = -\frac{1}{\xi}$$

Case-III

For $\sigma > 0$ then the Riccati equation can be written equation, considering integral constant as zero

$$\int d\xi = \int \frac{d\phi}{\sigma + \phi^2} \qquad \xi = \begin{cases} \frac{1}{\sqrt{\sigma}} \tan^{-1}\left(\frac{\phi}{\sqrt{\sigma}}\right) \\ -\frac{1}{\sqrt{\sigma}} \cot^{-1}\left(\frac{\phi}{\sqrt{\sigma}}\right) \end{cases}$$

Therefore, inverting above we get the following

$$\phi = \begin{cases} \sqrt{\sigma} \tan\left(\sqrt{\sigma}\,\xi\right) \\ -\sqrt{\sigma} \cot\left(\sqrt{\sigma}\,\xi\right) \end{cases}$$

The solution of the equation (2) will be expressed in-terms of $\phi$ in the form

$$u(\xi) = S = \sum_{i=0}^{n} a_i \phi^i \quad \text{.................(3)}$$

Where $n$ is obtained balancing the highest order derivative term and the non-linear terms in equation (2). Then put the value of $u(\xi)$ from (3) to the equation (2) and comparing the coefficient of $\phi^i$ the values of coefficients can be obtained.

Consider the equation $\frac{\partial^2 u}{\partial x^2} = u \frac{\partial u}{\partial t}$. Using the travelling wave transformation $\xi = kx + ct$ the considered equation reduces to

$$k^2 u_{\xi\xi} = c u u_\xi \quad \text{.................(4)}$$

Then using the transformation as defined in (3) in-terms of $\phi$ and using the chain rule of derivative we obtain

$$\frac{d}{d\xi} = \frac{d\phi}{d\xi}\frac{d}{d\phi} = \operatorname{sech}^2\xi \frac{d}{d\phi} = (1-\phi^2)\frac{d}{d\phi}$$

$$\frac{d^2}{d\xi^2} = \frac{d}{d\xi}\frac{d}{d\xi} = \frac{d}{d\phi}\left((1-\phi^2)\frac{d}{d\phi}\right)\frac{d\phi}{d\xi} = (1-\phi^2)^2\frac{d^2}{d\phi^2} - 2\phi(1-\phi^2)\frac{d}{d\phi}$$

Putting this transformation in equation (4) we get

$$k^2\left[(1-\phi^2)^2\frac{d^2S}{d\phi^2} - 2\phi(1-\phi^2)\frac{dS}{d\phi}\right] = cS(1-\phi^2)\frac{dS}{d\phi} \quad \ldots\ldots\ldots\ldots\ldots\ldots(5)$$

In the highest order derivative term the highest order of $\phi$ is $n+2$ and in the non-linear term it is $2n+1$. Equating them we get $n=1$. We have (3) as following

$$u(\xi) = S = a_0 + a_1\phi$$

Therefore

$$\frac{dS}{d\phi} = a_1 \qquad \text{and} \qquad \frac{d^2S}{d\phi^2} = 0$$

Putting the above value in equation (5) we get the following

$$k^2\left[2\phi(1-\phi^2)a_1\right] = c(a_0 + a_1\phi)(1-\phi^2)a_1$$

Comparing the like powers of $\phi$ from both sides of above expression we get

$$a_0 = 0 \qquad \text{and} \qquad a_1 = -\frac{2k^2}{c}$$

Hence the corresponding solution is

$$u(\xi) = S = -\frac{2k^2}{c}\phi = -\frac{2k^2}{c}\tanh(\xi)$$

that is

$$u(x,t) = -\frac{2k^2}{c}\tanh(kx + ct)$$

## 3.0 Fractional sub-equation method

The space and time fractional partial differential equation is of the form

$$L(u, u_t^\alpha, u_x^\alpha, u_x^{2\alpha}, \ldots) = 0 \qquad 0 < \alpha < 1 \ldots\ldots\ldots(6)$$

where $u = u(x,t)$ and $L$ is functions (linear or non-linear operator) of the $u$ and its derivatives. $\alpha$ is order of fractional derivative of Jumarie type, defined as follows

$$_0^J D_x^\alpha [f(x)] = f^{(\alpha)}(x) = \begin{cases} \dfrac{1}{\Gamma(-\alpha)} \int_0^x (x-\xi)^{-\alpha-1} f(\xi) d\xi, & \alpha < 0 \\ \dfrac{1}{\Gamma(1-\alpha)} \dfrac{d}{dx} \int_0^x (x-\xi)^{-\alpha} (f(\xi) - f(0)) d\xi, & 0 < \alpha < 1 \\ \left( f^{(\alpha-n)}(x) \right)^{(n)}, & n \le \alpha < n+1, \quad n \ge 1 \end{cases}$$

We consider that at $x < 0$ the function $f(x) = 0$ and also $f(x) - f(0) = 0$ for $x < 0$. The fractional derivative considered here in the fractional differential equation are using Jumarie [22] modified Riemann-Liouville (RL) derivative which is defined as above. The first expression above is fractional integration of Jumarie type. The modification by Jumarie is to carry RL fractional integration or RL fractional derivation is by forming a new function by offsetting the original function by subtraction the function value at the start point; and then operate the RL definition. The above defined fractional derivative has the following properties

$$_0^J D_x^\alpha [x^\gamma] = \frac{\Gamma(1+\gamma)}{\Gamma(1+\gamma-\alpha)} x^{\gamma-\alpha}, \qquad \gamma > 0,$$

$$_0^J D_x^\alpha [f(x) g(x)] = g(x) \left( _0^J D_x^\alpha f(x) \right) + f(x) \left( _0^J D_x^\alpha g(x) \right)$$

$$_0^J D_x^\alpha [f(g(x))] = f_g'(g(x)) \left( _0^J D_x^\alpha [g(x)] \right) = \left( _0^J D_g^\alpha [f(g(x))] \right) (g_x')^\alpha$$

The equation (7) is example of the fractional differential equation

$$\frac{\partial^{2\alpha} u}{\partial x^{2\alpha}} = u \frac{\partial^\alpha u}{\partial t^\alpha} \ldots\ldots\ldots(7)$$

Using the travelling wave transformation $\xi = kx + ct$, the equation (6) reduce to

$$L(u, u_\xi^\alpha, u_{\xi\xi}^{2\alpha}, \ldots) = 0 \qquad 0 < \alpha < 1 \ldots\ldots\ldots(8)$$

Where, $u = u(\xi)$ is function of one variable $\xi$

Using the same transformation, the equation (7) reduce to $\frac{d^{2\alpha}u}{d\xi^{2\alpha}} = u \frac{d^{\alpha}u}{d\xi^{\alpha}}$, as described below

Thus the operator $\frac{\partial^{\alpha}}{\partial t^{\alpha}}$, in terms of $\xi$ with $\xi = kx + ct$ is

$$\frac{\partial^{\alpha}}{\partial t^{\alpha}} \equiv \frac{\partial^{\alpha}}{\partial \xi^{\alpha}} \frac{\partial^{\alpha}\xi}{\partial t^{\alpha}} = c^{\alpha} \frac{d^{\alpha}}{d\xi^{\alpha}} = c\left(^{J}D_{\xi}^{\alpha}\right)$$

and similarly $\frac{\partial^{2\alpha}}{\partial x^{2\alpha}}$ is

$$\frac{\partial^{2\alpha}}{\partial x^{2\alpha}} \equiv \frac{\partial^{\alpha}}{\partial x^{\alpha}}\left(\frac{\partial^{\alpha}}{\partial x^{\alpha}}\right) = \frac{\partial^{\alpha}}{\partial \xi^{\alpha}} \frac{\partial^{\alpha}\xi}{\partial x^{\alpha}}\left(\frac{\partial^{\alpha}}{\partial \xi^{\alpha}} \frac{\partial^{\alpha}\xi}{\partial x^{\alpha}}\right) = \frac{d^{\alpha}}{d\xi^{\alpha}}k^{\alpha}\left(\frac{d^{\alpha}}{d\xi^{\alpha}}k^{\alpha}\right) = k^{2\alpha} \frac{d^{2\alpha}}{d\xi^{2\alpha}} = k^{2\alpha}\left(^{J}D_{\xi}^{2\alpha}\right)$$

The transformed equation is thus, $k^{2\alpha} \frac{\partial^{2\alpha}u}{\partial \xi^{2\alpha}} - c^{\alpha}u \frac{\partial^{\alpha}u}{\partial \xi^{\alpha}} = 0$ that is $L(u, u_{\xi}^{\alpha}, u_{\xi\xi}^{2\alpha}) = 0$.

The solution of the equation (8) is taken in the following series form

$$u(\xi) = \sum_{i=0}^{n} b_{i}\Phi^{i} \dots\dots\dots\dots\dots\dots\dots\dots\dots\dots\dots\dots\dots\dots(9)$$

Where $\Phi(\xi)$ satisfy the fractional Riccati equation

$$\Phi^{\alpha}(\xi) = \sigma + \Phi^{2}, \qquad 0 < \alpha < 1$$

Here $\sigma$ is a real constant. Say for $\alpha = 1$ we obtain the Riccati equation as $\frac{d\Phi}{d\xi} = \sigma + \Phi^{2}$. From this we write $\frac{d\Phi}{\sigma + \Phi^{2}} = d\xi$, integrating both sides once we get $\int \frac{d\Phi}{\sigma + \Phi^{2}} = \xi$. In terms of anti-derivative symbol, we write as $\xi = D_{\Phi}^{-1}\left[\frac{1}{\sigma + \Phi^{2}}\right]$. In similar way we write in terms of fractional integration (Jumarie type) the following

$$\xi = D_{\Phi}^{-\alpha}\left[\frac{1}{\sigma + \Phi^{2}}\right]$$

Similar to the previous method the coefficients can be obtained. Zhang et al [29] develop the method finding the solution fractional Riccati equation in the following form using generalized hyperbolic and trigonometric function

$$\Phi(\xi) = \begin{cases} -\sqrt{-\sigma}\tanh_\alpha\left(\sqrt{-\sigma}\xi\right) \\ -\sqrt{-\sigma}\coth_\alpha\left(\sqrt{-\sigma}\xi\right) \end{cases} \quad \text{for } \sigma < 0 \\ \begin{cases} \sqrt{\sigma}\tan_\alpha\left(\sqrt{\sigma}\xi\right) \\ -\sqrt{\sigma}\cot_\alpha\left(\sqrt{\sigma}\xi\right) \end{cases} \quad \text{for } \sigma > 0 \\ -\frac{\Gamma(1+\alpha)}{\xi^\alpha + \omega} \quad \text{for } \sigma = 0 \end{cases}$$

Where, the defined higher order transcendental functions are as follows

$$\tanh_\alpha(x) = \frac{\sinh_\alpha(x)}{\cosh_\alpha(x)} \qquad \coth_\alpha(x) = \frac{\cosh_\alpha(x)}{\sinh_\alpha(x)}$$

$$\sinh_\alpha(x) = \frac{E_\alpha(x^\alpha) - E_\alpha(x^\alpha)}{2} \qquad \cosh_\alpha(x) = \frac{E_\alpha(x^\alpha) + E_\alpha(x^\alpha)}{2}$$

$$\tan_\alpha(x) = \frac{\sin_\alpha(x)}{\cos_\alpha(x)} \qquad \cot_\alpha(x) = \frac{\cos_\alpha(x)}{\sin_\alpha(x)}$$

$$\sin_\alpha(x) = \frac{E_\alpha(ix^\alpha) - E_\alpha(ix^\alpha)}{2} \qquad \cos_\alpha(x) = \frac{E_\alpha(ix^\alpha) + E_\alpha(ix^\alpha)}{2}$$

Where $E_\alpha(z) = \sum_{k=0}^{\infty} \frac{z^\alpha}{\Gamma(1+k\alpha)}$ is the one parameter Mittag-Leffler function.

## 4.0 Nonlinear evolution equations (NEEs) of shallow water waves (SWW) partial and fractional differential equations

The nonlinear evolution equation of the shallow water wave equation is

$$u_{xt} + u_{xx} = u_{xxxy} + pu_x u_{xt} + qu_t u_{xx} \quad \text{...............................(10)}$$

Where $u = u(x, y, t)$. The Corresponding space and time fractional differential equation is

$$u_{xt}^{2\alpha} + u_{xx}^{2\alpha} = u_{xxxy}^{4\alpha} + pu_x^\alpha u_{xt}^{2\alpha} + qu_t^\alpha u_{xx}^{2\alpha} \quad \text{...........................(11)}$$

Where for $0 < \alpha < 1$, we have the following

$$u_{xt}^{2\alpha} = \frac{\partial^{2\alpha} u}{\partial t^\alpha \partial x^\alpha} \qquad u_{xxxy}^{4\alpha} = \frac{\partial^{4\alpha} u}{\partial x^{3\alpha} \partial y^\alpha} \qquad u_{xx}^{2\alpha} = \frac{\partial^{2\alpha} u}{\partial t^{2\alpha}}$$

## 5.0 Kadomtsev-Petviashvili (KP) partial and fractional differential equation

The KP equation was first written in 1970 by Soviet physicists Boris B. Kadomtsev and Vladimir I. Petviashvili. It came as a natural generalization of the KdV equation (derived by Korteweg and De Vries in 1895). Whereas in the KdV equation waves are strictly one-dimensional, in the KP equation this restriction is relaxed. The KP equation usually written as

$$u_{yy} = (u_t + 6uu_x + u_{xxx})_x \quad \quad (12)$$

Where $u = u(x, y, t)$. The Corresponding space and time fractional differential equation is

$$u_y^{2\alpha} = (u_t^{\alpha} + 6uu_x^{\alpha} + u_x^{3\alpha})_x^{\alpha} \quad \quad (13)$$

Where

$$u_y^{2\alpha} = \frac{\partial^{2\alpha} u}{\partial y^{2\alpha}} \quad \quad u_x^{3\alpha} = \frac{\partial^{3\alpha} u}{\partial t^{3\alpha}}, \quad 0 < \alpha < 1$$

## 6.0 Forth order Boussinesq equations partial and fractional differential equation

In fluid dynamics, the Boussinesq approximation for water waves is an approximation valid for weakly non-linear and fairly long waves. The approximation is named after Joseph Boussinesq, who first derived them in response to the observation by John Scott Russell of the wave of translation (also known as solitary wave or soliton). The 1872 paper of Boussinesq introduces the equations now known as the Boussinesq equations. The forth order Boussinesq equation can be written in the following form. Solitons are solitary waves with huge amplitude locally, and zero amplitude at far away, with very high travelling velocity

$$u_{tt} = u_{xx} + 3(u^2)_{xx} + u_{xxxx} \quad \quad (14)$$

Where $u = u(x,t)$. The Corresponding space and time fractional equation is

$$u_t^{2\alpha} = u_x^{2\alpha} + 3(u^2)_x^{2\alpha} + u_x^{4\alpha} \quad \quad (15)$$

Where

$$u_t^{2\alpha} = \frac{\partial^{2\alpha} u}{\partial t^{2\alpha}} \quad \quad u_x^{4\alpha} = \frac{\partial^{4\alpha} u}{\partial t^{4\alpha}}, \quad 0 < \alpha < 1$$

# 7.0 Solution of Nonlinear evolution equations (NEEs) of shallow water waves (SWW) using tanh method

Equation (10) that is $u_{xt} + u_{xx} = u_{xxxy} + pu_x u_{xt} + qu_t u_{xx}$ is non-linear equation in $u = u(x, y, t)$. Using travelling wave transformation in the form $\xi = ct + kx + my$ and following substitution

$$u_{xt} = cku_{\xi\xi} \qquad u_{xx} = k^2 u_{\xi\xi}$$
$$u_{xxxy} = k^3 m u_{\xi\xi\xi\xi} \qquad u_x u_{xt} = k^2 c u_\xi u_{\xi\xi}$$
$$u_t u_{xx} = k^2 c u_\xi u_{\xi\xi}$$

in equation (10), is reduced to the ordinary differential equation form

$$k^2 m u_{\xi\xi\xi\xi} + ck(p+q) u_\xi u_{\xi\xi} + (c+k) u_{\xi\xi} = 0 \quad \text{..............................(16)}$$

Since the solitary waves are localized in space therefore the solution and its derivatives at large distance from the pulse are known to be extremely small and they vanishes $\xi \to \pm\infty$ [28]. Integrating once w.r.t $\xi$ the above equation and by assuming the condition the derivative terms tend to zero as $\xi \to \infty$ we get

$$2k^2 m u_{\xi\xi\xi} + ck(p+q) u_\xi^2 + 2(c+k) u_\xi = 0 \quad \text{..............................(17)}$$

Here we consider the solution in the form
$$u(\xi) = S = a_0 \phi^0 + a_1 \phi + a_2 \phi^2 + \ldots + a_n \phi^n \quad \text{..............................(18)}$$

Where $\phi(\xi) = \tanh(\xi)$
Then

$$\frac{d}{d\xi} = (1 - \phi^2) \frac{d}{d\phi}$$

$$\frac{d^2}{d\xi^2} = (1 - \phi^2)^2 \frac{d^2}{d\phi^2} - 2\phi(1 - \phi^2) \frac{d}{d\phi}$$

$$\frac{d^3}{d\xi^3} = (1 - \phi^2)^3 \frac{d^3}{d\phi^3} - 6\phi(1 - \phi^2)^2 \frac{d^2}{d\phi^2} - 2\phi(1 - \phi^2)(1 - 3\phi^2) \frac{d}{d\phi}$$

$$\frac{d^4}{d\xi^4} = (1-\phi^2)^4 \frac{d^4}{d\phi^4} - 12\phi(1-\phi^2)^3 \frac{d^3}{d\phi^3} -$$

$$(36\phi^2 - 8)(1-\phi^2)^2 \frac{d^2}{d\phi^2} + 2\phi(1-\phi^2)(8-12\phi^2) \frac{d}{d\phi}$$

Putting the value of $u$ and using the transformation define above equation (17) reduce to the form

$$2k^2 m(1-\phi^2)^3 \frac{d^3 S}{d\phi^3} - 6\phi(1-\phi^2)^2 \frac{d^2 S}{d\phi^2}$$

$$-2(1-\phi^2)(1-3\phi^2)\frac{dS}{d\phi} + ck(p+q)\left[(1-\phi^2)\frac{dS}{d\phi}\right]^2$$

$$+2(c+k)(1-\phi^2)\frac{dS}{d\phi} = 0 \dots\dots\dots\dots\dots\dots\dots\dots\dots\dots\dots\dots\dots\dots\dots(19)$$

In the highest order derivative term the highest order of $\phi$ is $n+3$ and in the non-linear term it is $2n+2$. Equating them we get $n=1$.

We have (18) as $u(\xi) = S = a_0 + a_1\phi(\xi)$. Putting the value of $S$ in (19) we get

$$-4k^2 m(1-\phi^2)(1-3\phi^2)a_1 + ck(p+q)\left((1-\phi^2)a_1\right)^2 + 2(c+k)(1-\phi^2)a_1 = 0 \dots\dots\dots(20)$$

Comparing the powers of $\phi$ from both side the following relations obtained for $a_1 \neq 0$

$$\phi^0 : -2(c+k) + a_1 ck(p+q) = 4k^2 m$$
$$\phi^1 : 2(c+k) - 2a_1 ck(p+q) = -16k^2 m$$
$$\phi^3 : a_1 ck(\alpha + \beta) = 12k^2 m$$

Solving the above we get $a_1 = \frac{12km}{c(p+q)}$ and $(c+k) = 4k^2 m$.

Therefore the solution is

$$u(x, y, t) = a_0 + \frac{12km}{c(p+q)} \tanh\left(ct + kx + my\right)$$

for all $a_0 \in R$.

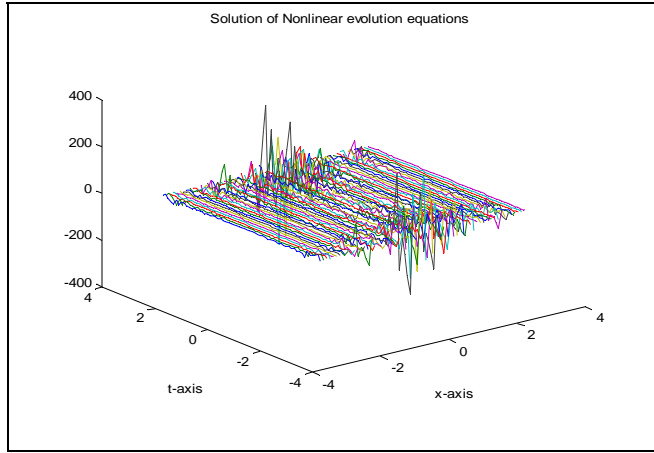

Figure 1: Graphical presentation of *u* for fixed *y* and considering $p = q = m = k = 1$, $c = 3$.

## 8.0 Solution of Nonlinear fractional evolution equations (NEEs) of shallow water waves (SWW) using Fractional Sub-equation method

Using the same travelling wave transformation as define in section-7.0 and with following substitutions

$$u_{xt}^{2\alpha} = c^\alpha k^\alpha u_\xi^{2\alpha}$$

$$u_{xx}^{2\alpha} = k^{2\alpha} u_\xi^{2\alpha}$$

$$u_{xxxy}^{4\alpha} = k^{3\alpha} m^\alpha u_\xi^{4\alpha}$$

$$u_x^\alpha = k^\alpha u_\xi^\alpha$$

$$u_t^\alpha = c^\alpha u_\xi^\alpha$$

in equation (16), that is $u_{xt}^{2\alpha} + u_{xx}^{2\alpha} = u_{xxxy}^{4\alpha} + p u_x^\alpha u_{xt}^{2\alpha} + q u_t^\alpha u_{xx}^{2\alpha}$ reduce to the form

$$k^{2\alpha} m^\alpha u_\xi^{4\alpha} + k^\alpha c^\alpha (p+q) u_\xi^{2\alpha} u_\xi^\alpha = (k^\alpha + c^\alpha) u_\xi^{2\alpha} \quad \text{...................(21)}$$

Consider the solution in the form (22)

$$u(\xi) = a_0 + a_1 \Phi + a_2 \Phi^2 + \ldots + a_n \Phi^n \quad \text{.......................................................(22)}$$

Where $\Phi$ is define in section 3.0 Comparing the highest powers $\Phi$ from the highest order derivative term and the non-linear derivative term from (21) and using *u* from (22) we get for $n + 4 = 3n + 3$ or $n = 1$. Therefore, $u = a_0 + a_1 \Phi(\xi)$. Putting the value of *u* from the above in (21) and comparing the like powers of $\Phi$ the following relations is obtained in simplest form and using the relation $a_1 \neq 0$

$$\Phi^1: \quad 2(c^\alpha + k^\alpha)\sigma = a_1 c^\alpha k^\alpha (p+q)\sigma^2 + 16k^{2\alpha} m^\alpha \sigma^2$$
$$\Phi^3: \quad 2(c^\alpha + k^\alpha) = 4a_1 c^\alpha k^\alpha (p+q)\sigma + 40k^{2\alpha} m^\alpha$$
$$\Phi^5: \quad 24k^{2\alpha} m^\alpha = 2a_1 c^\alpha k^\alpha (p+q)$$

Solving the above relation we get $a_1 = -\frac{12k^\alpha m^\alpha}{c^\alpha (\alpha+\beta)}$ and $(c^\alpha + k^\beta) + 4\sigma k^{2\alpha} m^\beta = 0$.
Therefore the solution is following

$$u(x,y,t) = \begin{cases} a_0 + \dfrac{12k^\alpha m^\alpha}{c^\alpha(p+q)}\sqrt{-\sigma}\,\tanh_\alpha\!\left(\sqrt{-\sigma}(ct+kx+my)\right) \\ a_0 + \dfrac{12k^\alpha m^\alpha}{c^\alpha(p+q)}\sqrt{-\sigma}\,\coth_\alpha\!\left(\sqrt{-\sigma}(ct+kx+my)\right) \end{cases} \text{for } \sigma < 0 \\ a_0 + \dfrac{12k^\alpha m^\alpha}{c^\alpha(p+q)}\dfrac{\Gamma(1+\alpha)}{(ct+kx+my)^\alpha + \omega} \quad \text{for } \sigma = 0, \quad \omega - \text{Cons tan t.} \\ \begin{cases} a_0 - \dfrac{12k^\alpha m^\alpha}{c^\alpha(p+q)}\sqrt{\sigma}\,\tan_\alpha\!\left(\sqrt{\sigma}(ct+kx+my)\right) \\ a_0 + \dfrac{12k^\alpha m^\alpha}{c^\alpha(p+q)}\sqrt{\sigma}\,\cot_\alpha\!\left(\sqrt{\sigma}(ct+kx+my)\right) \end{cases} \text{for } \sigma > 0 \end{cases} \quad \ldots\ldots(23)$$

For all $a_0 \in R$.

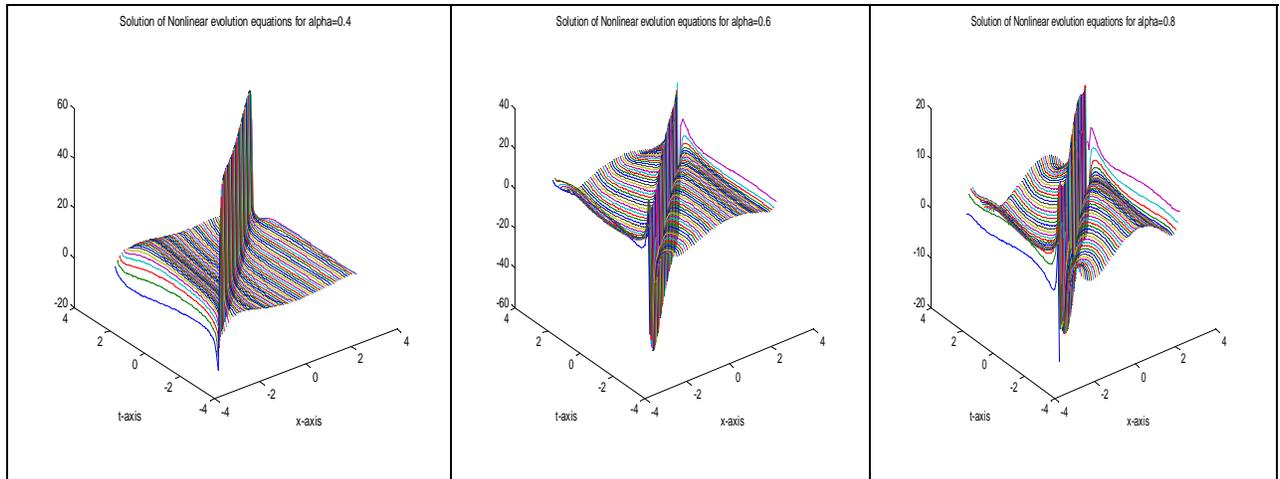

Figure 2: Graphical presentation of $u$ for fixed $y$ and considering $p = q = m = k = 1$, $c = 3$, for $\sigma = -1$.

# 9.0 Solution of Kadomtsev-Petviashvili (KP) equation using tanh method

Using the same travelling wave transformation as define in section-7.0 and $u_{yy} = m^2 u_{\xi\xi}$, $u_{xxx} = k^3 u_{\xi\xi\xi}$, $uu_x = kuu_\xi$ & $u_t = cu_\xi$ equation (12), that is $u_{yy} = (u_t + 6uu_x + u_{xxx})_x$ reduces to

$$m^2 u_{\xi\xi} = k(cu_\xi + 6kuu_\xi + k^3 u_{\xi\xi\xi})_\xi \quad \text{...................................................(24)}$$

Integrating both side two time w.r.t $\xi$ and then using the condition the solitary waves are localized in space that is $u, u_\xi, u_{\xi\xi}$..all tends to zero when $\zeta \to \infty$ [28]. We get after first integration $m^2 u_\xi = k(cu_\xi + 6kuu_\xi + k^3 u_{\xi\xi\xi})$ and after second integration we get

$$(m^2 - kc)u = 3k^2 u^2 + k^4 u_{\xi\xi} \quad \text{...................................................(25)}$$

We consider the solution as defined in equation (18), that is $u(\xi) = S = a_0 + a_1\phi + a_2\phi^2 + ... + a_n\phi^n$. Then putting the values of $u$ from (18) in (25) we get

$$(m^2 - kc)S = 3k^2 S^2 + k^4 \left( (1-\phi^2)^2 \frac{d^2 S}{d\phi^2} - 2\phi(1-\phi^2) \frac{dS}{d\phi} \right) \quad \text{...............................(26)}$$

Comparing the highest order derivative and the non-linear term we get for $n = 2$. Therefore, $u(\xi) = S = a_0 + a_1\phi(\xi) + a_2\phi^2(\xi)$, Putting the value of $u$ in (26) and comparing the like powers of $\phi$ the following relations are obtained

$$\phi^0: \quad (m^2 - kc)a_0 + 3k^2 a_0^2 + 2a_2 k^4 = 0$$
$$\phi^1: \quad (m^2 - kc)a_1 + 6k^2 a_0 a_1 - 2a_1 k^4 = 0$$
$$\phi^2: \quad (m^2 - kc)a_2 = 3k^2(a_1^2 + 2a_0 a_2) - 8a_2 k^4$$
$$\phi^3: \quad 6a_1 a_2 k^2 + 4a_1 k^4 = 0$$
$$\phi^4: \quad 3a_2^2 k^2 + 6a_2 k^4 = 0$$

Solving the above we get $a_0 = \frac{8k^4 + m^2 - kc}{6k^2}$, $a_1 = 0$ and $a_2 = -2k^2$ with the relation $9a_0^2 = 8k^2 + 2k^4$

Hence the solution is

$$u(x, y, t) = \frac{8k^4 + m^2 - kc}{6k^2} - 2k^2 \tanh^2(ct + kx + my)$$

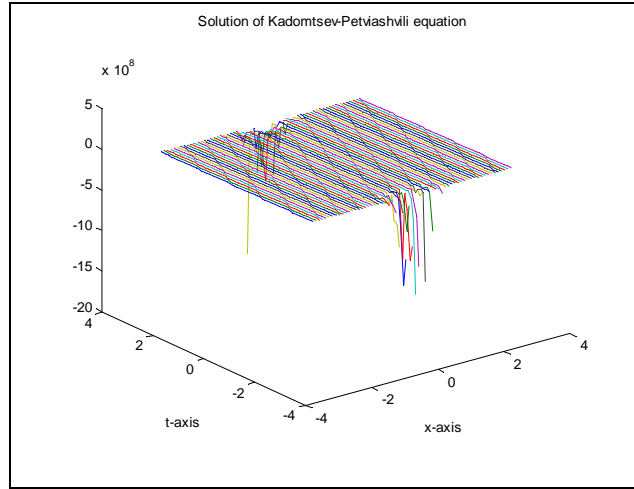

Figure 3: Graphical presentation of *u* for fixed *y* and considering $m = k = 1$, $c = 3.68$.

## 10.0 Solution of Kadomtsev-Petviashvili (KP) equation using Fractional Sub-equation method

The equation (13), that is $u_y^{2\alpha} = (u_t^{\alpha} + 6uu_x^{\alpha} + u_x^{3\alpha})_x^{\alpha}$ can be written in the following form

$$u_{tx}^{2\alpha} + 6(u_x^{\alpha})^2 + 6uu_{xx}^{2\alpha} + u_{xxxx}^{4\alpha} = u_{yy}^{2\alpha} \dots\dots\dots(27)$$

Using the same travelling wave transformation $\xi = ct + kx + my$ and with following substitutions

$$u_{tx}^{2\alpha} = kcu_{\xi\xi}^{2\alpha}$$
$$(u_x^{\alpha})^2 = k^2(u_{\xi}^{\alpha})^2$$
$$u_{xxxx}^{4\alpha} = k^{4\alpha}u_{\xi}^{4\alpha}$$
$$u_{yy}^{2\alpha} = m^2 u_{\xi}^{2\alpha}$$

equation (27) reduce to

$$(c^{\alpha}k^{\alpha} - m^{2\alpha})u_{\xi}^{2\alpha} + 6k^{2\alpha}(u_x^{\alpha})^2 + 6k^{2\alpha}uu_{\xi}^{2\alpha} + k^{4\alpha}u_{\xi}^{4\alpha} = 0\dots\dots\dots \quad (28)$$

We want to find the solution in the form as define in equation (22), that is $u(\xi) = a_0 + a_1\Phi + a_2\Phi^2 + \dots + a_n\Phi^n$, and then comparing the highest order derivative term and the non-linear term we get for $n = 2$. Therefore $u(\xi) = a_0 + a_1\Phi(\xi) + a_2\Phi^2(\xi)$. Putting in equation (27) and comparing the coefficients of powers of $\Phi(\xi)$ we get

$$\Phi^0: \quad 2(c^\alpha k^\alpha - m^{2\alpha})a_2\sigma^2 + 6k^{2\alpha}(2a_0 a_2 \sigma^2 + a_1^2 \sigma^2) + 16k^{4\alpha} a_2 \sigma^3 = 0$$

$$\Phi^1: \quad 2(c^\alpha k^\alpha - m^{2\alpha})a_1\sigma^2 + 6k^{2\alpha}(2a_0 a_1 \sigma + 6a_1 a_2 \sigma^2) + 16k^{4\alpha} a_1 \sigma^2 = 0$$

$$\Phi^2: \quad 8(c^\alpha k^\alpha - m^{2\alpha})a_2\sigma + 6k^{2\alpha}(8a_0 a_2 \sigma + 4a_1^2 \sigma + 2a_2^2 \sigma) + 136k^{4\alpha} a_2 \sigma^2 = 0$$

$$\Phi^3: \quad 2(c^\alpha k^\alpha - m^{2\alpha})2a_1\sigma + 6k^{2\alpha}(2a_0 a_1 + 18a_1 a_2 \sigma) + 40k^{4\alpha} a_1 = 0$$

$$\Phi^4: \quad 6(c^\alpha k^\alpha - m^{2\alpha})a_2 + 6k^{2\alpha}(6a_0 a_2 + 16a_2^2 \sigma + 3a_1^2) + 240k^{4\alpha} a_2 \sigma = 0$$

$$\Phi^5: \quad 72k^{2\alpha} a_1 a_2 + 24 a_1 k^{4\alpha} = 0$$

$$\Phi^6: \quad 60k^{2\alpha} a_2^2 + 120 a_2 k^{4\alpha} = 0$$

Solving we get $a_0 = \frac{m^{2\alpha} - 8\sigma k^{4\alpha} - k^\alpha c^\alpha}{6k^2}$, $a_1 = 0$ and $a_2 = -2k^{2\alpha}$. Hence the solution is following

$$u(\xi) = \begin{cases} \dfrac{m^{2\alpha} - 8\sigma k^{4\alpha} - k^\alpha c^\alpha}{6k^2} + 2k^{2\alpha} \sigma \tanh_\alpha^2 \left(\sqrt{-\sigma}(ct + kx + my)\right) \\ \dfrac{m^{2\alpha} - 8\sigma k^{4\alpha} - k^\alpha c^\alpha}{6k^2} + 2k^{2\alpha} \sigma \coth_\alpha^2 \left(\sqrt{-\sigma}(ct + kx + my)\right) \end{cases} \quad \text{for} \quad \sigma < 0$$

$$\dfrac{m^{2\alpha} - 8\sigma k^{4\alpha} - k^\alpha c^\alpha}{6k^2} - 2k^{2\alpha} \left( \dfrac{(\Gamma(1+\alpha))^2}{\left((ct + kx + my)^\alpha + \omega\right)^2} \right) \quad \text{for} \quad \sigma = 0, \omega - \text{Cons tan t.}$$

$$\begin{cases} \dfrac{m^{2\alpha} - 8\sigma k^{4\alpha} - k^\alpha c^\alpha}{6k^2} - 2k^{2\alpha} \sigma \tan_\alpha^2 \left(\sqrt{\sigma}(ct + kx + my)\right) \\ \dfrac{m^{2\alpha} - 8\sigma k^{4\alpha} - k^\alpha c^\alpha}{6k^2} - 2k^{2\alpha} \sigma \cot_\alpha^2 \left(\sqrt{\sigma}(ct + kx + my)\right) \end{cases} \quad \text{for} \quad \sigma > 0$$

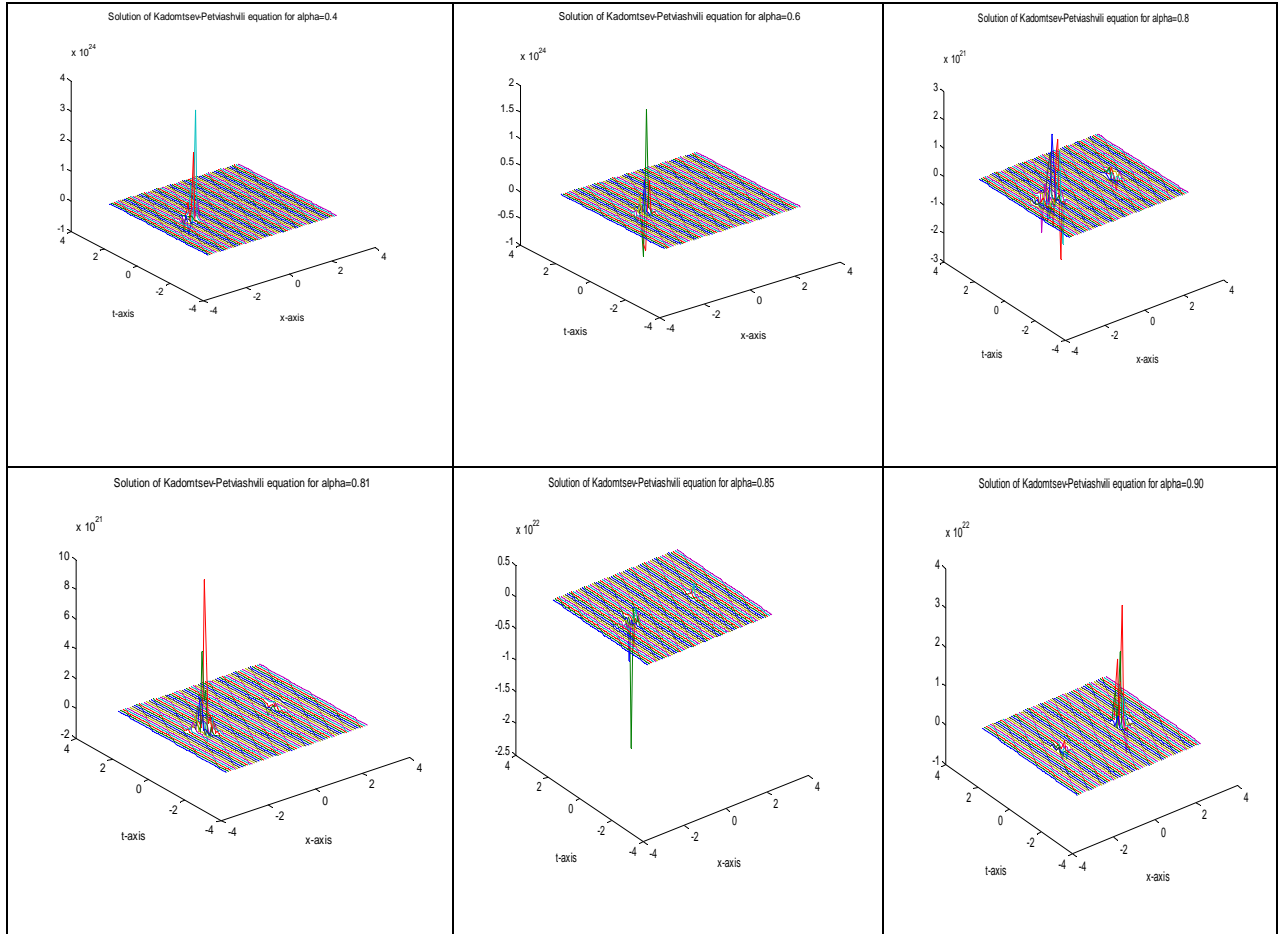

Figure 4: Graphical presentation of *u* for fixed *y* and considering *m* = *k* = 1, *c* = 3.68, for $\sigma = -1$.

## 11.0 Solution of Forth order Boussinesq equation using tanh Method

Equation (14), that is non-linear equation in $u = u(x,t)$. Using travelling wave transformation in the form $\xi = kx + ct$ and with following substitution

$$u_{tt} = c^2 u_{\xi\xi}$$
$$u_{xx} = k^2 u_{\xi\xi}$$
$$(u^2)_x = k^2 (u^2)_\xi$$
$$u_{xxxx} = k^4 u_{\xi\xi\xi\xi}$$

in equation (14) the equation reduced to the form

$$(c^2 - k^2)u_{\xi\xi} = 3k^2 (u^2)_{\xi\xi} + k^4 u_{\xi\xi\xi\xi} \ldots\ldots\ldots\ldots\ldots\ldots\ldots\ldots(29)$$

Here we consider the solution as defined in equation (18), that is $u(\xi) = S = a_0 + a_1\phi + a_2\phi^2 + ... + a_n\phi^n$. Putting the values of $u$ from (9) in equation (19) and using the transformation as define in section 7.0 we get

$$(c^2 - k^2)\left[(1-\phi^2)^2 \frac{d^2S}{d\phi^2} - 2\phi(1-\phi^2)\frac{dS}{d\phi}\right] =$$

$$3k^2\left[\{2(1-\phi^2)^2 \frac{d^2S}{d\phi^2} - 4\phi(1-\phi^2)\frac{dS}{d\phi}\}S\right] + 2(1-\phi^2)^2\left(\frac{dS}{d\phi}\right)^2 +$$

$$k^4 \begin{pmatrix} (1-\phi^2)^4 \frac{d^4S}{d\phi^4} - 12\phi(1-\phi^2)^3 \frac{d^3S}{d\phi^3} \\ +(1-\phi^2)^2(36\phi^2 - 8)\frac{d^2S}{d\phi^2} + 2\phi(1-\phi^2)(8-12\phi^2)\frac{dS}{d\phi} \end{pmatrix}$$

$$= 0 ................................(30)$$

Comparing the highest order derivative term and the non-linear term we obtain $n = 2$. Therefore $u(\xi) = S = a_0 + a_1\phi(\xi) + a_2\phi^2(\xi)$. Putting this value of $S$ in equation (30) and comparing the powers of $\phi$ we get

$$\phi^0: \quad 2(c^2 - k^2)a_2 = 3k^2(4a_0a_2 + 2a_1^2) - 16a_2k^4$$

$$\phi^1: \quad -2(c^2 - k^2)a_1 = 3k^2(-4a_0a_1 + 12a_1a_2) + 16a_1k^4$$

$$\phi^2: \quad -8(c^2 - k^2)a_2 = 3k^2(-16a_0a_2 + 4a_1^2) + 136a_2k^4$$

$$\phi^3: \quad 2(c^2 - k^2)a_1 = 3k^2(4a_0a_1 - 36a_1a_2) - 40a_1k^4$$

$$\phi^4: \quad 6(c^2 - k^2)a_2 = 3k^2(-26a_1^2 + 12a_0a_2) - 240a_2k^4$$

$$\phi^5: \quad 24a_1k^4 + 72a_1a_2k^2 = 0$$

$$\phi^6: \quad 60a_2k^2(a_2 + 2k^2) = 0$$

Solving the above set of equation we get $a_0 = \frac{8k^4 + c^2 - k^2}{6k^4}$, $a_1 = 0$ and $a_2 = -2k^2$. Hence the solution is

$$u(x, y, t) = \frac{8k^4 + c^2 - k^2}{6k^4} - 2k^2 \tanh^2(ct + kx)$$

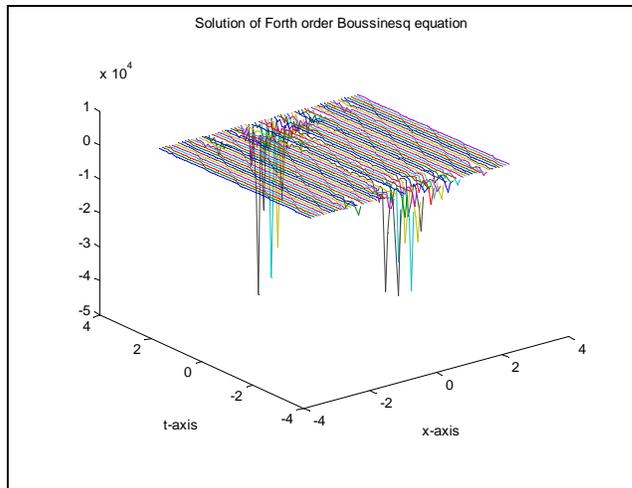

Figure 5: Graphical presentation of *u* for fixed *y* and considering c = k = 1.

## 12.0: Solution of Forth order fractional Boussinesq equation using Fractional Sub-equation method

Using the same type of travelling wave transformation $\xi = ct + kx$ and with following substitution

$$u_t^{2\alpha} = c^{2\alpha} u_\xi^{2\alpha},$$
$$u_x^{2\alpha} = k^{2\alpha} u_\xi^{2\alpha},$$
$$(u^2)_x^{2\alpha} = k^{2\alpha} (u^2)_\xi^{2\alpha}$$
$$u_x^{4\alpha} = k^{4\alpha} u_\xi^{4\alpha}$$

the fractional Boussinesq equation (15) which is $u_t^{2\alpha} = u_x^{2\alpha} + 3(u^2)_x^{2\alpha} + u_x^{4\alpha}$ reduce to the following form

$$(c^{2\alpha} - k^{2\alpha}) u_\xi^{2\alpha} = 3k^{2\alpha} (u^2)_\xi^{2\alpha} + k^{4\alpha} u_\xi^{4\alpha} \quad \ldots\ldots\ldots\ldots\ldots\ldots(31)$$

We want to find the solution in the form as define in equation (22), comparing the highest order derivative term and the non-linear term we get $n = 2$. Therefore, $u(\xi) = a_0 + a_1 \Phi(\xi) + a_2 \Phi^2(\xi)$. Comparing the like powers of $\Phi(\xi)$ we get the following $\Phi^0:$  $\quad 2(c^{2\alpha}-k^{2\alpha})a_2\sigma^2=3k^{2\alpha}\sigma^2(4a_0a_2+2a_1^2)+16a_2k^{4\alpha}\sigma^3$

$\Phi^1:$  $\quad 2\sigma(c^{2\alpha}-k^{2\alpha})a_1=12k^{2\alpha}(a_0a_1\sigma+3a_1a_2\sigma^2)+16a_1\sigma^2k^{4\alpha}$

$\Phi^2:$  $\quad 8\sigma(c^{2\alpha}-k^{2\alpha})a_2=24k^{2\alpha}(2a_0a_2+a_1^2)\sigma+36\sigma^2a_2^2k^{2\alpha}+136a_2k^{4\alpha}\sigma^2$

$\Phi^3:$  $\quad 2(c^{2\alpha}-k^{2\alpha})a_1=120k^{2\alpha}a_0a_1+108a_1a_2\sigma k^{2\alpha}+40a_1\sigma k^{4\alpha}$

$\Phi^4:$  $\quad 6(c^{2\alpha}-k^{2\alpha})a_2=18k^{2\alpha}(a_1^2+2a_0a_2)+96k^{2\alpha}a_2^2\sigma+240\sigma a_2k^{4\alpha}$

$\Phi^5:$  $\quad 24a_1k^{4\alpha}+72a_1a_2k^{2\alpha}=0$

$\Phi^6:$  $\quad 60a_2k^{2\alpha}(a_2+2k^{2\alpha})=0$

Solving the above system of equation we get $a_0=\frac{c^{2\alpha}-8\sigma k^{4\alpha}-k^{2\alpha}}{6k^2}$, $a_1=0$ and $a_2=-2k^{2\alpha}$. Hence the solution is

$$u(x,y)=\begin{cases} \dfrac{c^{2\alpha}-8\sigma k^{4\alpha}-k^{2\alpha}}{6k^2}+2k^{2\alpha}\sigma\tanh_\alpha^2\left(\sqrt{-\sigma}(ct+kx)\right) \\ \dfrac{c^{2\alpha}-8\sigma k^{4\alpha}-k^{2\alpha}}{6k^2}+2k^{2\alpha}\sigma\coth_\alpha^2\left(\sqrt{-\sigma}(ct+kx)\right) \end{cases} \text{for } \sigma<0 \\ \dfrac{c^{2\alpha}-8\sigma k^{4\alpha}-k^{2\alpha}}{6k^2}-2k^{2\alpha}\left(\dfrac{(\Gamma(1+\alpha))^2}{\left((ct+kx)^\alpha+\omega\right)^2}\right) \quad \text{for } \sigma=0, \quad \omega-\text{Cons}\tan\text{t.} \\ \begin{cases} \dfrac{c^{2\alpha}-8\sigma k^{4\alpha}-k^{2\alpha}}{6k^2}-2k^{2\alpha}\sigma\tan_\alpha^2\left(\sqrt{\sigma}(ct+kx)\right) \\ \dfrac{c^{2\alpha}-8\sigma k^{4\alpha}-k^{2\alpha}}{6k^2}-2k^{2\alpha}\sigma\cot_\alpha^2\left(\sqrt{\sigma}(ct+kx)\right) \end{cases} \text{for } \sigma>0$$

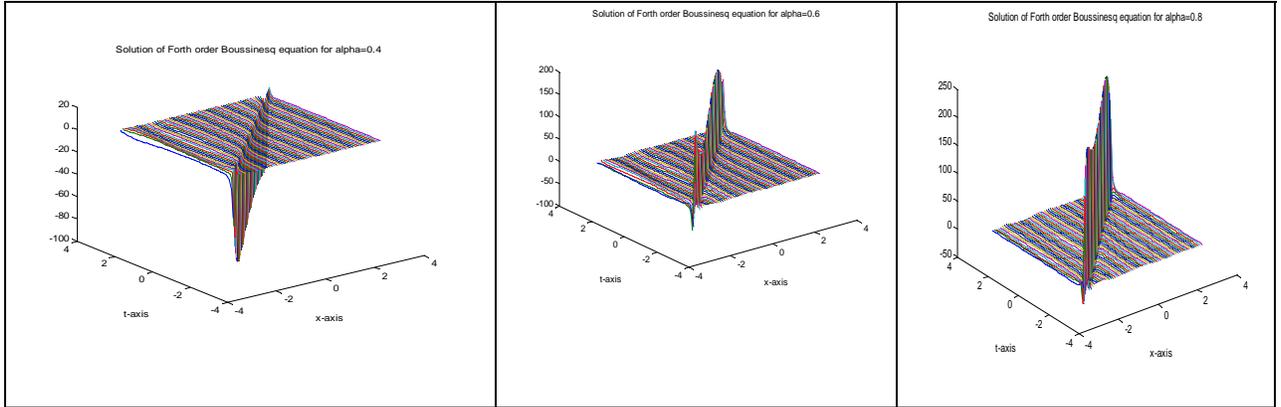

Figure 6: Graphical presentation of $u$ for fixed $y$ and considering $c=k=1$, for $\sigma=-1$.

From figure 1-6 it is clear that solutions of the fractional differential equations depend on order of fractional derivative.

# 13.0 Conclusion

In all the three cases the solution obtained by tanh-method coincides with one solution of obtained by Fractional Sub-Equation method, when fractional order tends to one. Thus the solution of the non-linear partial differential equation and the non-linear fractional differential equation is obtained in analytic form. However, the solution obtained by fractional sub-equation method is more general compare to the tanh-method, as it encompasses the fractional generalization of the non-linear partial differential equations.